\documentclass{elsart}

\usepackage{amsmath,amssymb,amsopn,amsfonts,mathrsfs,amsbsy,amscd}

\usepackage[english,francais]{babel}

\newcommand{\esp}{\quad\mbox{and}\quad}

\newcommand{\la}{\lambda}

\newcommand{\ad}{{\mathrm{ad}}}

\newtheorem{theorem}{Theorem}[section]

\newtheorem{e-proposition}[theorem]{Proposition}
\newtheorem{corollary}[theorem]{Corollary}
\newtheorem{e-definition}[theorem]{Definition\rm}

\newcommand{\G}{{\mathfrak{g}}}
\newcommand{\prs}{\langle\;,\;\rangle}
\newtheorem{theoreme}{Th\'eor\`eme}[section]

\newtheorem{proposition}[theoreme]{Proposition}
\newtheorem{corollaire}[theoreme]{Corollaire}

\def\og{\leavevmode\raise.3ex\hbox{$\scriptscriptstyle\langle\!\langle$~}}
\def\fg{\leavevmode\raise.3ex\hbox{~$\!\scriptscriptstyle\,\rangle\!\rangle$}}

\journal{arxiv}
\begin{document}
\centerline{}
\begin{frontmatter}


\selectlanguage{english}
\title{Lorentzian Flat Lie Groups Admitting a Timelike Left-Invariant Killing Vector Field.}


\selectlanguage{english}
\author[authorlabel1]{Hicham Lebzioui}
\ead{hlebzioui@gmail.com}

\address[authorlabel1]{Facult\'{e} des sciences de Mekn\`{e}s. Morocco}


\medskip

\begin{abstract}
\selectlanguage{english}
We call a connected Lie group endowed with a left-invariant Lorentzian flat metric
Lorentzian flat Lie group. In this Note, we determine all
Lorentzian flat Lie groups admitting a timelike left-invariant Killing
vector field. We show that  these Lie groups
are 2-solvable and unimodular and hence geodesically complete.
Moreover, we show that a Lorentzian flat Lie group $(\mathrm{G},\mu)$ admits
a timelike left-invariant Killing vector field if and only if $\mathrm{G}$
admits a left-invariant Riemannian metric which has the same Levi-Civita
connection of $\mu$. Finally, we give an useful characterization of left-invariant
pseudo-Riemannian flat metrics on Lie groups $\mathrm{G}$ satisfying the property: for any couple of left invariant vector fields $X$ and $Y$ their Lie bracket
 $[X,Y]$ is a linear combination of $X$ and $Y$.

\vskip 0.5\baselineskip

\selectlanguage{francais}
\noindent{\bf R\'esum\'e} \vskip 0.5\baselineskip \noindent {\bf
Groupes de Lie Lorentziens Plats Admettant Un Champ de Killing
Invariant \`{a} Gauche de Type Temps. }\\ On appelle groupe de Lie
lorentzien plat un groupe de Lie connexe muni d'une m\'{e}trique lorentzienne plate invariante \`{a} gauche. Dans cette Note, on d\'{e}termine tout les groupes de Lie lorentziens plats qui admettent un champ de Killing invariant \`{a} gauche de type temps. On montre que de tels groupes de Lie  sont 2-r\'{e}solubles et unimodulaires et donc g\'{e}od\'{e}siquement complets. On montre aussi qu'un groupe de Lie lorentzien plat $(\mathrm{G},\mu)$ admet un champ de Killing invariant \`{a} gauche de type temps si et seulement si $\mathrm{G}$ poss\'{e}de une m\'{e}trique riemannienne invariante \`{a} gauche qui a la m\^{e}me connexion de Levi-Civita de $\mu$. Finalement, on donne une caract\'{e}risation simple des m\'{e}triques pseudo-riemanniennes plates invariantes \`{a} gauche sur les groupes de Lie $\mathrm{G}$ qui v\'{e}rifient la propri\'et\'e suivant: pour tout couple de champ de vecteurs invariant \`a gauche $X$ et $Y$, leur crochet de Lie ${[X,Y]}$ est combinaison lin\'{e}aire de ${X}$ et ${Y}$.

\end{abstract}
\end{frontmatter}

\selectlanguage{francais}
\section*{Version fran\c{c}aise abr\'eg\'ee}
Un groupe de Lie pseudo-riemannien   est un groupe de Lie $\mathrm{G}$ muni d'une m\'{e}trique pseudo-riemannienne  invariante \`{a} gauche $\mu$.  Si la m\'{e}trique est riemannienne (resp. de signature $(-,+,\dots,+)$), le groupe de Lie est dit riemannien (resp. lorentzien). Soit $(\mathrm{G},\mu)$ un groupe de Lie pseudo-riemannien  et
  $(\G,\prs)$ son alg\`{e}bre de Lie o\`{u} $\prs$ est la valeur de $\mu$ en l'\'{e}l\'{e}ment neutre $e$. On d\'{e}note par $S(\G)$ la sous-alg\`{e}bre de $\G$ d\'{e}finie par $S(\G)=\{u\in \G/\ad_u+\ad_u^*=0\}$ o\`{u} $\ad_u^*$ est l'adjoint de $\ad_u$ relativement \`{a} $\prs$. Il est facile de voir qu'un champ de vecteur $X$ sur $\mathrm{G}$ invariant \`{a} gauche est de Killing si et seulement si $X(e)\in S(\G)$.
  La connexion de Levi-Civita de $\mu$ d\'{e}finie sur $\G$ un produit dit de Levi-Civita  par la relation \eqref{lv}. Il est connu qu'un groupe de Lie pseudo-riemannien plat est g\'{e}od\'{e}siquement complet si et seulement si le groupe est unimodulaire \cite{AM}.  Dans \cite{M}, J. Milnor a caract\'{e}ris\'{e} les groupes de Lie riemanniens plats, et tout ces groupes sont g\'{e}od\'{e}siquement complets. Dans \cite{BB} Theorem 3.1, les auteurs ont donn\'{e} une  formulation plus pr\'{e}cise du th\'{e}or\`{e}me de Milnor. Dans le cas lorentzien, les seuls groupes de Lie lorentziens plats connus sont ceux nilpotents  \cite{AM}, et ceux dont le centre de leur alg\`ebre de Lie est d\'{e}g\'{e}n\'{e}r\'{e} \cite{ABL}. Tous ces groupes de Lie lorentziens sont g\'{e}od\'{e}siquement complets. Rappelons qu'un champ de vecteur $X$  dans une vari\'{e}t\'{e} lorentzienne $(M,\mu)$ est dit de type temps si, pour tout $p\in M$, $\mu_p(X_p,X_p)<0$. Plusieurs auteurs se sont int\'{e}r\'{e}ss\'{e}s aux vari\'{e}t\'{e}s lorentziennes qui admettent un champ de Killing de type temps. Cette condition  donne des informations importantes sur la structure de la vari\'{e}t\'{e} (voire par exemple \cite{Sa}). Le premier objectif de cette Note est de d\'{e}montrer le th\'{e}or\`{e}me suivant:
\begin{theoreme}Soit $(\mathrm{G},\mu)$ un groupe de Lie lorentzien. Alors $\mathrm{G}$ est plat et admet un champ de Killing invariant
\`{a} gauche de type temps  si et seulement si son alg\`{e}bre de
Lie $\G$ se d\'{e}compose en une somme orthogonale $\G=S(\G)\oplus [\G,\G]$,
avec $[\G,\G]$ ab\'{e}lien, $S(\G)$ ab\'{e}lienne et contienne un vecteur de type temps.
Dans ce cas $[\G,\G]$ est de dimension paire, et le produit de Levi-Civita
est donn\'{e} par \eqref{eql}.
\end{theoreme}
On d\'{e}duit de ce th\'{e}or\`{e}me les corollaires suivants:
\begin{corollaire}
Tout groupe de Lie lorentzien plat qui admet un champ de Killing invariant \`{a} gauche de type temps est 2-r\'{e}soluble et unimodulaire et donc il est g\'{e}od\'{e}siquement complet.
\end{corollaire}
\begin{corollaire}
Soit $(\mathrm{G},\mu)$ un groupe de Lie lorentzien plat. Alors $\mathrm{G}$ admet un champ de Killing invariant \`{a} gauche de type temps si et seulement si $\mathrm{G}$ admet une m\'{e}trique riemannienne invariante \`{a} gauche qui a la m\^{e}me connexion de Levi-Civita que $\mu$.
\end{corollaire}
On d\'{e}note par $\mathcal{C}$ la classe des groupes de Lie $\mathrm{G}$ non ab\'{e}lien et v\'{e}rifiant, pour tout $x,y\in\G$,  le crochet $[x,y]$ est combinaison lin\'{e}aire de $x$ et $y$.
Le deuxi\`{e}me objectif de cette Note est de donner une condition n\'{e}cessaire et suffisante pour qu'une m\'{e}trique pseudo-riemannienne invariante \`{a} gauche sur un groupe de Lie dans $\mathcal{C}$ soit plate.
\begin{theoreme}
Soit $\mathrm{G}\in \mathcal{C}$ et $\mu$ une m\'{e}trique pseudo-riemannienne invariante \`{a} gauche sur $\mathrm{G}$. Alors  $\mu$ est plate si et seulement si la restriction de $\prs$ \`{a} l'id\'{e}al d\'{e}riv\'{e} $[\G,\G]$ est d\'{e}g\'{e}n\'{e}r\'{e}.
\end{theoreme}
On d\'{e}duit ais\'{e}ment le corollaire suivant.\begin{corollaire}
Soit $\mathrm{G}\in \mathcal{C}$. Alors  toute m\'{e}trique pseudo-riemannienne invariante \`{a} gauche $\mu$ sur $\mathrm{G}$ dont la restriction de $\prs$ \`{a} $[\G,\G]$ est d\'{e}g\'{e}n\'{e}r\'{e} est incompl\`{e}te.\end{corollaire}
\selectlanguage{english}
\section{Introduction and main results}
\label{}
A Lie group $\mathrm{G}$ together with a left-invariant pseudo-Riemannian
metric is called \emph {pseudo-Riemannian Lie group}. When the metric is definite positive or of signature $(-,+,\ldots,+)$ the group is called Riemannian or Lorentzian.  Let $(\mathrm{G},\mu)$ be a pseudo-Riemannian Lie group and $\G$ its Lie algebra endowed with the inner product $\prs=\mu(e)$.
 The Levi-Civita connection of
$(\mathrm{G},\mu)$ defines a product $(u,v)\longrightarrow uv$ on $\G$ called
\emph{Levi-Civita product} given by:
\begin{equation}\label{lv}2\langle
uv,w\rangle=\langle[u,v],w\rangle-\langle[v,w],u\rangle+\langle[w,u],v\rangle.
\end{equation}
For any $u\in \G$, we denote by $\mathrm{L}_u$ and $\mathrm{R}_u$ respectively the left
multiplication and the right multiplication on $\G$ given by
$\mathrm{L}_u(v)=uv$ and $\mathrm{R}_u(v)=vu$.  For any $u\in \G$, $\mathrm{L}_u$ is
skew-symmetric with respect to $\prs$, and $\ad_u=\mathrm{L}_u-\mathrm{R}_u$
where $\ad_u:\G\longrightarrow \G$ is given by $\ad_u(v)=[u,v]$.\\
The curvature of $(\G,\prs)$ is given by
$K(u,v)=\mathrm{L}_{[u,v]}-[\mathrm{L}_u,\mathrm{L}_v].$
If $K$ vanishes identically  $(\mathrm{G},\mu)$ is called
\emph{pseudo-Riemannian flat Lie group}. Let $S(\G)$ be the Lie subalgebra of $\G$ defined by
$S(\G)=\{u\in \G/\ad_u+\ad_u^*=0\},$
where $\ad_u^*$ is the adjoint of $\ad_u$ with respect to
$\prs$. It is easy to see that a  left invariant vector field $X$ on $\mathrm{G}$ is a Killing vector field iff $X(e)\in S(\G)$.\\
It is a well-known result that a pseudo-Riemannian flat Lie group is geodesically complete if and only if it is unimodular \cite{AM}.
In \cite{M}, J. Milnor characterized Riemannian flat Lie groups. In
 \cite{ABL,BB} there is a more precise  version  of
Milnor's Theorem: A Lie group is a Riemannian flat Lie group if and only if
its Lie algebra $\G$ splits as an orthogonal direct sum:
$\G=S(\G)\oplus [\G,\G]$, where $S(\G)$ and $[\G,\G]$ are abelian.
Moreover, in this case the dimension of $[\G,\G]$ is even and the
Levi-Civita product is given by:\begin{equation}\label{eql}
\mathrm{L}_a=\left\{\begin{array}{ccc}\ad_a&if&a\in S(\G),\\0&if&a\in
[\G,\G].\end{array}\right.\end{equation}
In the Lorentzian case,  the only known Lorentzian flat Lie groups are
those nilpotent \cite{AM} and those with degenerate center
\cite{ABL}. Their Lie algebras are determined by the double extension process, and
the metric in both cases is complete.

Recall that  a vector field $X$ on a Lorentzian manifold $(M,g)$ is called timelike (resp. spacelike, null) if $g_p(X_p,X_p)<0$ (resp. $g_p(X_p,X_p)>0$, $g_p(X_p,X_p)=0$) for any $p\in M$.
Many authors in mathematics and physics are interested in the
existence of timelike Killing vector field in a Lorentzian manifold.
This condition implies interesting information on the structure of
the manifold. In General Relativity, if this condition is satisfied,
then the Lorentzian manifold is called stationary, and it is known
that any compact stationary manifold is geodesically complete
\cite{Sa}.

Our first and main purpose in this paper is to determine
Lorentzian flat Lie groups admitting a timelike left-invariant
killing vector field. More precisely, we will prove the following result.
\begin{theorem}\label{th1}Let $(\mathrm{G},\mu)$ be a Lorentzian Lie group. Then it is flat and carries a timelike left invariant killing
vector field  if and only if its Lie algebra $\G$ splits as an
orthogonal direct sum $\G=S(\G)\oplus [\G,\G]$, where $[\G,\G]$ is
abelian, and $S(\G)$ is abelian and contains a timelike vector.
Moreover, in this case the dimension of $[\G,\G]$ is even and the
Levi-Civita product is given by \eqref{eql}.
\end{theorem}
\begin{corollary}\label{co1}
A Lorentzian flat Lie group which admits a timelike
left-invariant killing vector field is 2-solvable and unimodular and
hence geodesically complete.
\end{corollary}
If the Lie group is non abelian, then we find a family of non
compact stationary manifolds that are geodesically complete.
\begin{corollary}\label{co2}
Let $(\mathrm{G},\mu)$ be a Lorentzian flat Lie group.
The existence of timelike left-invariant killing vector field on $\mathrm{G}$
is equivalent to the existence of left-invariant Riemannian metric
on $\mathrm{G}$ with the same Levi-Civita connection of $\mu$.
\end{corollary}On the other hand, we define the class $\mathcal{C}$ consisting of those Lie groups $G$ which are non commutative and their Lie algebras satisfy: for all $x$ and $y$ in $\mathrm{G}$,  $[x,y]$ is a linear combination
of $x$ and $y$.
This class $\mathcal{C}$  has
been studied by several authors.  J. Milnor
in \cite{M} showed that every left-invariant Riemannian metric on
$\mathrm{G}\in \mathcal{C}$ has constant negative sectional curvature. K. Nomizu in \cite{N}
showed that every left invariant Lorentzian metric on $\mathrm{G}\in \mathcal{C}$ has
constant sectional curvature which can be positive, negative or
null.\\
The following theorem gives an useful characterization of flat pseudo-Riemannian left invariant metric on a Lie group belonging to $\mathcal{C}$.
\begin{theorem}\label{th2}Let $\mathrm{G}$ be a non commutative Lie group that belongs to
$\mathcal{C}$, and $\mu$ a left-invariant pseudo-Riemannian metric
on $\mathrm{G}$. Then
$\mu$ is flat if and only if the
restriction of $\mu(e)$ to $[\G,\G]$ is degenerate.
\end{theorem}We deduce the following corollary.\begin{corollary}Let $\mathrm{G}$ in
$\mathcal{C}$. Then $\mathrm{G}$ is non unimodular and hence
any left-invariant pseudo-Riemannian metric $\mu$ on $\mathrm{G}$ such that the restriction of $\mu(e)$ to $[\G,\G]$ is degenerate is geodesically incomplete.\end{corollary}
\section{Proofs of Theorems}
This section is devoted to the proofs of the results of the last section. The following proposition proved in \cite{ABL} Lemma 3.1 will be useful later.

\begin{proposition}\label{pr}
Let $(\mathrm{G},\mu)$ be a Riemannian or Lorentzian flat Lie group, then:
$$S(\G)=(\G\G)^\bot=\{u\in \G/\mathrm{R}_u=0\}.$$In particular $S(\G)$ is abelian.
\end{proposition}

\subsection{Proof of Theorem \ref{th1}}

Let $(\G,\prs)$ be a flat Lorentzian Lie algebra. Suppose that there exists $u\in S(\G)=(\G\G)^\perp$ such that $\langle u,u\rangle<0$. Then the restriction of $\prs$ to $u^\perp$ is definite positive and so for $S(\G)^\perp=\G\G\subset u^\perp$. Then
$ \G=S(\G)\oplus \G\G$ and the restriction of $\prs$ to $S(\G)$ is nondegenerate Lorentzian.
 Let us show that $\mathfrak{h}=\G\G$ is abelian and $\mathfrak{h}=[\G,\G]$. It is clear that  $\mathfrak{h}$ is a Riemannian flat Lie algebra and hence
 \[ \mathfrak{h}=S(\mathfrak{h})\oplus[\mathfrak{h},\mathfrak{h}], \]where $S(\mathfrak{h})$ and $[\mathfrak{h},\mathfrak{h}]$ are abelian. Moreover, there exist $\la_1,\ldots,\la_p\in S(\mathfrak{h})^*\setminus\{0\}$ and orthonormal basis $(e_1,\ldots,e_{2p})$ of $[\mathfrak{h},\mathfrak{h}]$ such that
 \[ \forall i=1,\ldots,p,\forall s\in S(\mathfrak{h}),\quad  \ad_{s}(e_{2i-1})=\la_i(s)e_{2i} \esp
 \ad_{s}(e_{2i})=-\la_i(s)e_{2i-1}.\]
 On the other hand, for any $s\in S(\G)$, $\ad_s$ which is skew-symmetric leaves $\mathfrak{h}$ invariant and so it leaves  $[\mathfrak{h},\mathfrak{h}]$ also invariant and hence $S(\mathfrak{h})$. Moreover, if $x\in S(\mathfrak{h})$ satisfies $[x,s]$ for any $s\in S(\G)$, then $\ad_x$ is skew-symmetric and then $x\in S(\G)$ and hence $x=0$. So there exists also $\mu_1,\ldots,\mu_q\in S(\G)^*\setminus\{0\}$ and orthonormal basis $(f_1,\ldots,f_{2q})$ of $S(\mathfrak{h})$ such that
  \[ \forall i=1,\ldots,q,\forall s\in S(\G),\quad  \ad_{s}(f_{2i-1})=\mu_i(s)f_{2i} \esp
  \ad_{s}(f_{2i})=-\mu_i(s)f_{2i-1}.\]
Suppose that $\mathfrak{h}$ is non abelian. Without loss of generality we can suppose that $[f_1,e_1]=\la_1(f_1)e_2\not=0$. We have also $[f_1,e_2]=-\la_1(f_1)e_1\not=0$. Moreover, there exists $s\in S(\G)$ such that $[s,f_1]=\mu_1(s)f_2\not=0$. We have necessarily $[f_2,e_1]=\la_1(f_2)e_2\not=0$. Otherwise, if $[f_2,e_1]=0$ then $[f_2,e_2]=0$ and hence
\begin{eqnarray*}
\ad_{[s,f_2]}(e_1)&=&-\mu_1(s)\la_1(f_1)e_2\\&=&
[\ad_{s},\ad_{f_2}](e_1)\\
&=&[[s,e_1],f_2]\\
&=&\sum_{i=3}^{2p}\alpha_ie_i.
\end{eqnarray*}
This is impossible since $\mu_1(s)\la_1(f_1)\not=0$. On the other hand, we have
\begin{eqnarray*}
0&=&[e_1,[f_1,s]]+[f_1,[s,e_1]]+[s,[e_1,f_1]]\\
&=&-\mu_1(s)[e_1,f_2]+[f_1,\sum_{i=2}^{2p}\beta_i e_i]-\la_1(f_1)[s,e_2]\\
&=&\mu_1(s)\la_1(f_2)e_2+\gamma_1 e_1+\sum_{i=3}^{2p}\gamma_i e_i.
\end{eqnarray*}This is impossible and hence $\mathfrak{h}$ is abelian.\\
Now it is obvious that $[\G,\G]\subset\mathfrak{h}$ and $[\mathfrak{h},\mathfrak{h}]\subset[\G,\G]$. On the other hand, for any $i=1,\ldots,2q$, there exists $s\in S(\G)$ such that $[s,f_{2i-1}]=\mu_i(s)f_{2i}\not=0$ and $[s,f_{2i}]=-\mu_i(s)f_{2i-1}\not=0$. This shows that $f_{2i-1},f_{2i}\in[\G,\G]$.
 Conversely, if the Lie algebra $\G$ of the Lorentzian Lie group
$(\mathrm{G},\mu)$ splits as an orthogonal direct sum $\G=S(\G)\oplus [\G,\G]$
where $S(\G)$ and $[\G,\G]$ are abelian, and there exists a timelike
vector in $S(\G)$ then the Levi-Civita product is given by \eqref{eql}
 which implies that $(\mathrm{G},\mu)$ is a
Lorentzian flat Lie group which admits a timelike left-invariant
Killing vector field.\hfill$\square$

\subsection{Proof of Corollary \ref{co2}}

Let $(\mathrm{G},\mu)$ be a Lorentzian flat Lie group and
$(\G,\prs)$ its Lorentzian flat Lie algebra.
If $(\mathrm{G},\mu)$ admits a timelike left-invariant Killing vector field,
then according to Theorem \ref{th1}, $\G=S(\G)\oplus[\G,\G]$ where $S(\G)$ is nondegenerate Lorentzian and $[\G,\G]$ is nondegenerate Euclidean. Choose any definite positive product $\prs_s$ on $S(\G)$ and define $\prs'$ on $\G$ as the orthogonal sum of $\prs_s$ and the restriction of $\prs$ to $[\G,\G]$. It is easy to check that the left invariant Riemannian metric on $\mathrm{G}$ associated to $\prs'$ has the same Levi-Civita connection as $\mu$.\\
Conversely, suppose that $\mathrm{G}$ possesses a left-invariant Riemannian metric
$\mu'$ with the same Levi-Civita connection of $\mu$. Put $\prs=\mu(e)$ and $\prs'=\mu'(e)$. Since $\mu$ and $\mu'$ have the same Levi-Civita connection then, according to Proposition \ref{pr} $\mu$ and $\mu'$ have also the same space of left invariant Killing vector fields identified to $S(\G)$. Then, by applying both \cite{ABL} Theorem 3.1 and Theorem \ref{th1}, we get $\G=S(\G)\oplus[\G,\G]$ and this splitting is orthogonal wit respect to $\prs$ and $\prs'$. If the restriction of $\prs$ to $S(\G)$ is Lorentzian then $\mu$ admits obviously a left invariant timelike Killing vector field. Suppose now that the restriction of $\prs$ to $[\G,\G]$ is Lorentzian. Then there exists a basis $(e_1,\ldots,e_{p})$ of $[\G,\G]$ which is orthogonal with respect to $\prs$ and orthonormal with respect to $\prs'$ and such that, $\langle e_{i},e_{i}\rangle>0$ for $i=1,\ldots,p-1$ and  $\langle e_{p},e_{p}\rangle<0$. Let $s\in S(\G)$. We have $\langle[s,e_{p}],e_{p}\rangle=0$ and for any $i=1,\ldots,p-1$,
\begin{eqnarray*}
\langle[s,e_{p}],e_i\rangle&=&\langle[s,e_{p}],e_i\rangle'\langle e_i,e_i\rangle\\
&=&-\langle[s,e_{i}],e_p\rangle\\
&=&-\langle[s,e_{i}],e_p\rangle'\langle e_p,e_p\rangle\\
&=&\langle[s,e_{p}],e_i\rangle'\langle e_p,e_p\rangle.
\end{eqnarray*}Since $\langle e_p,e_p\rangle<0$ and $\langle e_i,e_i\rangle>0$ then
$\langle[s,e_{p}],e_i\rangle=0$ and hence $[s,e_p]=0$ for any $s\in S(\G)$ and hence $e_p$ lies in the center of $\G$ and thus $e_p\in S(\G)$ which is impossible. This completes the proof of the corollary.\hfill$\square$

\subsection{Proof of Theorem \ref{th2}}

Let $(\mathrm{G},\mu)$ be a pseudo-Riemannian Lie group and
$(\G,\prs)$ its pseudo-Riemannian Lie algebra.\\
Milnor in \cite{M} has shown that $\mathrm{G}\in \mathcal{C}$ if and only if
$\G$ contains an abelian ideal $\mathcal{U}$ of codimension 1 and an
element $b\notin \mathcal{U}$ such that $[b,x]=x$ for every $x\in
\mathcal{U}$. Note that $\mathcal{U}=[\G,\G]$.\\
If $\mu$ is flat, then since $\mathrm{G}$ is non unimodular, thus according
to \cite{ABL} Proposition 3.1, we deduce that the restriction of $\prs$
to $[\G,\G]$ is degenerate.\\
Conversely, if the restriction of $\prs$ to $[\G,\G]$ is
degenerate, then there exists $e\in [\G,\G]$ such that $\langle
e,x\rangle=0$ for any $x\in [\G,\G]$. Let $y\in \G$ such that $\langle
e,y\rangle\neq 0$. We put:
$$d=\frac{y}{\langle y,e\rangle}-\frac{1}{2}\frac{\langle y,y\rangle}{\langle y,e\rangle^2}e,$$ then
$\langle d,e\rangle=1$ and $\langle d,d\rangle=0$.\\
The restriction of $\prs$ to $\mathrm{span}\{e,d\}$ is non
degenerate, thus if we put $\mathcal{B}=\mathrm{span}\{e,d\}^\bot$ then we
have $\G=\mathbb{R}e\oplus\mathcal{B}\oplus \mathbb{R}d$ where
$[\G,\G]=\mathbb{R}e\oplus \mathcal{B}$.\\
We have $d=\alpha b+u_0$ where $\alpha\in\mathbb{R}\setminus\{0\}$ and $u_0\in
[\G,\G]$. Then the only non-vanishing brackets are $[d,e]=\alpha e$, $[d,u]=\alpha u$ for any $u\in [\G,\G].$
Using the equation \ref{lv}, we deduce that the Levi-Civita product is given by:
$$L_e=0\ ,\ de=\alpha e\ ,\ dd=-\alpha d\ ,\ du=ue=0\ ,\ ud=-\alpha
u\ ,\ uu'=\alpha \langle u,u'\rangle e$$ for
$u,u'\in\mathcal{B}$.\\
From these relations we obtain $K(e,d)=K(e,u)=K(d,u)=K(u,u')=0$ for $u,u'\in\mathcal{B}$, which
implies that the curvature $K$ vanishes identically.\hfill$\square$




\section*{Acknowledgements}
I would like to give my sincere thanks to Professor Mohamed Boucetta for his helpful discussions and encouragements, and Professor Malika Ait Ben Haddou for her support and guidance.

\end{document}